\documentclass[a4paper]{article}
\usepackage{amsmath}\usepackage{epsf,amsfonts,amsthm}
\usepackage{fontenc,indentfirst, delarray,amsfonts,amsmath,amssymb}

\newcommand{\be}{\begin{equation}}
\newcommand{\ee}{\end{equation}}
\newcommand{\bea}{\begin{eqnarray*}}
\newcommand{\eea}{\end{eqnarray*}}
\newcommand{\w}{\wedge}
\newcommand{\p}{\partial}

\newcommand{\raa}{\rightarrow}
\newcommand{\Ci}{C^{\infty}}
\newcommand{\E}{\ell}
\newcommand{\N}{{\mathbf N}}
\newcommand{\Z}{{\mathbf Z}}
\newcommand{\Q}{{\mathbf Q}}
\newcommand{\R}{{\mathbf R}}

\newcommand{\D}{(x_1^2+x_2^2)x_3}

\newcommand{\lp}{\left(}
\newcommand{\rp}{\right)}

\newcommand{\sh}{\sharp}
\newcommand{\op}[1]{\!\!\mathop{\rm ~#1}\nolimits}

\mathchardef\za="710B  
\mathchardef\zb="710C  
\mathchardef\zg="710D  
\mathchardef\zd="710E  
\mathchardef\zve="710F 
\mathchardef\zz="7110  
\mathchardef\zh="7111  
\mathchardef\zy="7112 
\mathchardef\zi="7113  
\mathchardef\zk="7114  
\mathchardef\zl="7115  
\mathchardef\zm="7116  
\mathchardef\zn="7117  
\mathchardef\zx="7118  
\mathchardef\zp="7119  
\mathchardef\zr="711A  
\mathchardef\zs="711B  
\mathchardef\zt="711C  
\mathchardef\zu="711D  
\mathchardef\zf="711E 
\mathchardef\zq="711F  
\mathchardef\zc="7120  
\mathchardef\zw="7121  
\mathchardef\ze="7122  
\mathchardef\zvy="7123  
\mathchardef\zvw="7124  
\mathchardef\zvr="7125 
\mathchardef\zvs="7126 
\mathchardef\zvf="7127  
\mathchardef\zG="7000  
\mathchardef\zD="7001  
\mathchardef\zY="7002  
\mathchardef\zL="7003  
\mathchardef\zX="7004  
\mathchardef\zP="7005  
\mathchardef\zS="7006  
\mathchardef\zU="7007  
\mathchardef\zF="7008  
\mathchardef\zW="700A  

\begin{document}

\title{On a general approach to the \\formal cohomology of quadratic Poisson
structures\footnote{The research of N. Poncin was supported by
grant R1F105L10}}
\author{Mohsen Masmoudi\footnote{Institut Elie Cartan, Universit\'e Henri Poincar\'e, B.P. 239, F-54 506 Vandoeuvre-les-Nancy Cedex, France, E-mail: Mohsen.Masmoudi@iecn.u-nancy.fr},
Norbert Poncin\footnote{University of Luxembourg, Campus
Limpertsberg, Mathematics Laboratory, 162A, avenue de la
Fa\"iencerie, L-1511 Luxembourg City, Grand-Duchy of Luxembourg,
E-mail: norbert.poncin@uni.lu}}\maketitle

\newtheorem{re}{Remark}
\newtheorem{theo}{Theorem}
\newtheorem{prop}{Proposition}
\newtheorem{lem}{Lemma}
\newtheorem{cor}{Corollary}
\newtheorem{ex}{Example}

\begin{abstract} We propose a general approach to the formal Poisson cohomology
of $r$-matrix induced quadratic structures, we apply this device
to compute the cohomology of structure $2$ of the Dufour-Haraki
classification, and provide complete results also for the
cohomology of structure $7$.\\

{\bf Key-words:} Poisson cohomology, formal cochain, quadratic
Poisson tensor, $r$-matrix

{\bf 2000 Mathematics Subject Classification:} 17B63, 17B56
\end{abstract}

\section{Introduction}

The Poisson cohomology was defined by Lichnerowicz in \cite{AL}.
It plays an important role in Poisson Geometry, gives several
information on the geometry of the manifold, is closely related to
the classification of singularities of Poisson structures,
naturally appears in infinitesimal deformations of Poisson tensors
[second cohomology group (non-trivial infinitesimal deformations),
third cohomology group (obstructions to extending a first order
deformation to a formal deformation)], and is exploited in the
study and the classification of star-products.

The case of regular Poisson manifolds is discussed in \cite{PX}
and \cite{IV}, some upshots regarding the Poisson cohomology of
Poisson-Lie groups can be found in \cite{GW}. Many results have
been proven only in two dimensions. Nakanishi \cite{NN} has
computed, using an earlier idea of Vaisman, the cohomology for
plane quadratic structures. But also more recent papers by Monnier
\cite{PM2} and Roger and Vanhaecke \cite{RV} are confined to
dimension $2$. Ginzburg \cite{VLG} has studied a spectral
sequence, Poisson analog of the Leray spectral sequence of a
fibration, which converges to the Poisson cohomology of the
manifold. Related cohomologies, the Nambu-Poisson cohomology
(which generalizes in a certain sense the Poisson cohomology in
dimension 2) \cite{PM1} or the tangential Poisson cohomology
(which governs tangential star products and is involved in the
Poisson cohomology) \cite{AnGa} have (partially) been computed.
Roytenberg has computed the cohomology on the 2-sphere for special
covariant structures \cite{DR} and Pichereau has taken an interest
in Poisson (co)homology and isolated singularities \cite{AP}.

In Deformation Quantization we are interested in the formal
Poisson cohomology, where "formal" means that cochains are
multi-vector fields with coefficients in formal series. Let us
emphasize that the formal Poisson cohomology associated to a
Poisson manifold ($M,P$), where the Poisson tensor $P$ gives
Kontsevich's star-product $*_K$, is linked with the Hochschild
cohomology of the associative algebra ($\Ci(M)[[\hbar]],*_K$) of
formal series in $\hbar$ with coefficients in $\Ci(M)$.

In \cite{PM3}, Monnier has computed the formal cohomology of
diagonal Poisson structures.

The aim of this work is to provide a general approach to the
formal Poisson cohomology of a broad set of isomorphism classes of
quadratic structures and to illustrate this modus operandi through
its application to two of the most demanding classes of
the Dufour-Haraki classification (DHC, see \cite{DH}).\\

Our procedure applies to the quadratic Poisson tensors $\zL$ that
read as linear combinations of wedge products of mutually
commuting linear vector fields. In the three-dimensional Euclidean
setting this means that $$\zL=aY_2\w Y_3+bY_3\w Y_1+c Y_1\w
Y_2\quad a,b,c\in\R,$$ where the vector fields $Y_1,Y_2,Y_3$ have
linear coefficients and meet the condition $[Y_i,Y_j]=0$. Let us
stress that most of the structures of the DHC have this {\it
admissible} form. The advantage of the just defined family of
admissible tensors is readily understood. If we substitute the
$Y_i$ for the standard basic vector fields $\p_i=\p/\p_{x_i}$, the
cochains assume---broadly speaking---the shape $\sum f\mathbf{Y}$,
where $f$ is a function and $\mathbf{Y}$ is a wedge product of
basic fields $Y_i$. Then the Lichnerowicz-Poisson (LP) coboundary
operator $[\zL,\cdot]$ is just
\[[\zL,f\mathbf{Y}]=[\zL,f]\w\mathbf{Y}.\] More precisely, \[[\zL,f]=\sum_iX_i(f)Y_i,\]
where $X_1=aY_2-cY_3, X_2=bY_3-aY_1,X_3=cY_1-bY_2.$ So the
coboundary of a $0$-cochain is a kind of gradient, and, as easily
checked, the coboundaries of a $1$- and a $2$-cochain, decomposed
in the $Y_i$-basis are a sort of curl and of divergence
respectively. As in the above ``gradient'' the operators $X_i$ act
in these ``curl'' and ``divergence'' as substitutes for the usual
partial derivatives. The preceding simplification of the
Lichnerowicz-Poisson coboundary operator is of course not
restricted to the three-dimensional context.

Let us emphasize that, if cochains are decomposed in the new
$Y_i$-induced bases, their coefficients are rational with fixed
denominator. Hence a natural injection of the {\it real cochain
space} ${\cal R}$ in a larger {\it potential cochain space} ${\cal
P}$. As a matter of fact, such a potential cochain is implemented
by a real cochain if and only if specific divisibility conditions
are satisfied. This observation directly leads to a {\it
supplementary cochain space} ${\cal S}$ of ${\cal R}$ in ${\cal
P}$. It is possible to heave space ${\cal S}$ into the category of
differential spaces. Since ${\cal P}$ is as ${\cal R}$ a complex
for the LP coboundary operator, we end up with a short exact
sequence of differential spaces and an exact triangle in
cohomology. It turns out that ${\cal S}$-cohomology and ${\cal
P}$-cohomology are less intricate than ${\cal R}$-cohomology, but
are important stages on the way to ${\cal R}$-cohomology, i.e. to
the cohomology of the
considered admissible quadratic Poisson structure.\\

In this work we provide explicit results---obtained by the just
depicted method---for the cohomologies of class 2 and class 7 of
the DHC. Observe that the representative
$$\zL_7=b(x_1^2+x_2^2)\p_1\w\p_2+((2b+c)x_1-ax_2)x_3\p_2\w\p_3+(ax_1+(2b+c)x_2)x_3\p_3\w\p_1$$
of class 7 reduces to the representative of class 2 for parameter
value $c=0$. Actually computations are quite similar in both
cases, hence we refrained from publishing those pertaining to case 7.\\

We now detail the Lichnerowicz-Poisson cohomology of structure
class 2. Remark first that $$\Lambda_2=2bY_2\w Y_3+aY_3\w
Y_1+bY_1\w Y_2,$$ where $Y_1=x_1\p_1+x_2\p_2$,
$Y_2=x_1\p_2-x_2\p_1$, and $Y_3=x_3\p_3$. As case $b=0$ has been
studied in \cite{PM3}, we assume that $b\neq 0$. We denote the
determinant $(x_1^2+x_2^2)x_3$ of the basic fields $Y_i$ by $D$.
The results of this article will entail that the abundance of
cocycles that do not bound is tightly related with closeness of
the considered Poisson tensor to Koszul-exactness. If $a=0$,
structure $\zL_2$ is exact and induced by $bD$. The algebra of
Casimir functions is generated by $1$ if $a\neq 0$ and by $D$ if
$a=0$. When rewording this statement and writing $\zL$ instead of
$\zL_2$, we get
$$H^0(\zL)=\op{Cas}(\zL)=\left\{\begin{array}{ll}\R,&\mbox{if}\quad a\neq 0,\\\bigoplus_{m=0}^{\infty}\R D^m,&
\mbox{if}\quad a=0.\end{array}\right.$$ Remind now that our
Poisson tensor is built with linear infinitesimal Poisson
automorphisms $Y_i$. It follows that the wedge products of the
$Y_i$ constitute ``a priori'' privileged cocycles. Of course,
2-cocycle $\zL_2$ itself, is a linear combination of such
privileged cocycles. Moreover, the curl or modular vector field
reads here $K(\zL)=a(2Y_3-Y_1)$ and is thus also a combination of
privileged cocycles. As the LP cohomology is an associative graded
commutative algebra, the first cohomology group of $\zL_2$ is easy
to conjecture:
$$H^1(\zL)=\bigoplus_{i}\op{Cas}(\zL)Y_i.$$ It is well-known that the singularities
of the investigated Poisson structure appear in the second and
third cohomology spaces. Observe that the singular points of
structure $\zL_2$ are the annihilators of $D'=x_1^2+x_2^2$. Note
also that any homogeneous polynomial $P\in\R[x_1,x_2,x_3]$ of
order $m$ reads
$$P=\sum_{\E=0}^m x_3^{\E}P_{\E}(x_1,x_2)=\sum_{\E=0}^m x_3^{\E}\lp
D'\cdot Q_{\E}+A_{\E}x_1^{m-\E}+B_{\E}x_1^{m-\E-1}x_2\rp,$$ with
self-explaining notations. Hence, the algebra of ``polynomials''
of the affine algebraic variety of singularities is
$$\bigoplus_{m=0}^{\infty}\bigoplus_{\E=0}^{m} x_3^{\E}\lp \R x_1^{m-\E}+\R
x_1^{m-\E-1}x_2\rp.$$ The third cohomology group contains a part
of this formal series. More precisely,
$$H^3(\zL)=\op{Cas}(\zL)\;Y_{123}\oplus\left\{\begin{array}{ll}\R\;\p_{123},&\mbox{if}\quad
a\neq 0,\\\bigoplus_{m=0}^{\infty}\R\;
x_3^m\p_{123}\oplus\bigoplus_{m=0}^{\infty}x_1^m\lp\R x_1+\R
x_2\rp\p_{123},& \mbox{if}\quad a=0,\end{array}\right.$$ where
$Y_{123}$ (resp. $\p_{123}$) means $Y_1\w Y_2\w Y_3$ (resp.
$\p_1\w\p_2\w\p_3$). The reader might object that the
``mother-structure'' $\zL$ is symmetric in $x_1,x_2$ and that
there should therefore exist a symmetric ``twin-cocycle''
$\bigoplus_{m=0}^{\infty}x_2^{m}(\R x_1+\R x_2)\p_{123}$. This
cocycle actually exists, but is---as easily checked---cohomologous
to the visible representative. Finally, the second cohomology
space reads
$$\begin{array}{ll}H^2(\zL)=&\op{Cas}(\zL)Y_{23}\oplus\op{Cas}(\zL)Y_{31}\oplus\op{Cas}(\zL)Y_{12}\\
&\oplus\left\{\begin{array}{ll}/\;,&\mbox{if}\quad a\neq
0,\\\bigoplus_{m=0}^{\infty}\R\;
x_3^m\p_{12}\oplus\bigoplus_{m=0}^{\infty}x_1^{m-1}\lp\R
x_1\p_{23}+\R (x_1\p_{31}+mx_2\p_{23})\rp,& \mbox{if}\quad
a=0.\end{array}\right.\end{array}$$ For $m\ge 1$, the last cocycle
has the form $$\lp \R x_1^m +\R
x_1^{m-1}x_2\rp\p_{23}+\lp\int\p_{x_2}\!\!\lp \R x_1^m +\R
x_1^{m-1}x_2\rp dx_1\rp\p_{31}$$ and is thus also induced by
singularities.\\

Our paper is organized as follows. Sections $2$-$8$ deal with the
explicit computation of the cohomology of tensor $\zL_2$. In
Section \ref{coho7} we completely describe the cohomology of
structure $\zL_7$. Finally, the last section is devoted to
explanations regarding the connection of the Poisson structures
that are accessible to our method with $r$-matrix induced Poisson
structures. Upshots relating to the admissible quadratic tensors
different from $\zL_2$ and $\zL_7$ are being published separately.

\section{Simplified differential}


In this and the following sections, we compute the formal
cohomology of the second quadratic Poisson structure of the
Dufour-Haraki classification. This structure reads
\[\Lambda=b(x_1^2+x_2^2)\p_1\w\p_2+x_3(2bx_1-ax_2)\p_2\w\p_3+x_3(ax_1+2bx_2)\p_3\w\p_1,\]
where $b\neq 0$ (if $b=0$, we recover the case studied in
\cite{PM3}).

As mentioned above, we get \be\Lambda=bY_1\w Y_2+2bY_2\w
Y_3+aY_3\w Y_1\label{La},\ee if we set \be Y_1=x_1\p_1+x_2\p_2,
Y_2=x_1\p_2-x_2\p_1, Y_3=x_3\p_3.\label{Y}\ee We also know that it
is interesting---in order to simplify the coboundary operator---to
write the cochains in terms of these fields. Let us recall that in
the formal setting, the cochains are the multi-vector fields with
coefficients in $\R[[x]]$, the space of formal series in
$x=(x_1,x_2,x_3)$ with coefficients in $\R$.

Note first that if $D=\D$ and if $x_{ij}=x_ix_j,$
\[\p_1=\frac{1}{D}(x_{13}Y_1-x_{23}Y_2),
\p_2=\frac{1}{D}(x_{23}Y_1+x_{13}Y_2),
\p_3=\frac{1}{D}(x_1^2+x_2^2)Y_3.\] An arbitrary $0$-cochain
$C^0=\sum_{I\in\N^3}c_Ix^{I}$ ($I=(i_1,i_2,i_3), c_I\in\R,
x^{I}=x^{i_1,i_2,i_3}=x_1^{i_1}x_2^{i_2}x_3^{i_3}$) can be written
\[
C^0=\sum_{J\in\N^3}\frac{x^J}{D}\alpha_J=:\frac{\zs}{D}=:\varsigma,\]
with self-explaining notations. Similarly, $1$-, $2$-, and
$3$-cochains $C^1=\zs_1\p_1+\zs_2\p_2+\zs_3\p_3,$
$C^2=\zs_1\p_{23}+\zs_2\p_{31}+\zs_3\p_{12},$ $C^3=\zs \p_{123},$
where $\zs_i,\zs\in\R[[x]]$ ($i\in\{1,2,3\}$) and where for
instance $\p_{23}=\p_2\w\p_3$, read in terms of the $Y_i$
($i\in\{1,2,3\}$),
$C^1=\varsigma_1Y_1+\varsigma_2Y_2+\varsigma_3Y_3,$
$C^2=\varsigma_1Y_{23}+\varsigma_2Y_{31}+\varsigma_3Y_{12},$ and
$C^3=\varsigma Y_{123}$, again with obvious notations. If
$[\cdot,\cdot]$ denotes the Schouten-bracket, the coboundary of
$C^0$ is given by

\[[\zL,C^0]=X_1(\varsigma)Y_1+X_2(\varsigma)Y_2+X_3(\varsigma)Y_3,\] where
\[X_1=bY_2-aY_3,X_2=-bY_1+2bY_3,X_3=aY_1-2bY_2.\]
Set $\nabla = \sum_i X_i(.)Y_i$. A short computation then shows
that
\[[\zL,C^0]=\nabla C^0,[\zL,C^1]=\nabla\w
C^1,[\zL,C^2]=\nabla.C^2,\mbox{ and }[\zL,C^3]=0,\] where the
r.h.s. have to be viewed as notations that give the coefficients
of the coboundaries in the $Y_i$-basis. For instance,
$[\zL,C^2]=(\sum_iX_i(\varsigma_i))Y_{123}$. We easily verifie
that
\[X_1\lp \frac{x^J}{D}\rp =[b\lp
j_2\frac{x_1}{x_2}-j_1\frac{x_2}{x_1}\rp-a(j_3-1)]\frac{x^J}{D},\]
\[X_2(\frac{x^J}{D})=b[2j_3-(j_1+j_2)]\frac{x^J}{D},\] and
\[X_3(\frac{x^J}{D})=[a(j_1+j_2-2)-2b\lp
j_2\frac{x_1}{x_2}-j_1\frac{x_2}{x_1}\rp]\frac{x^J}{D}.\]

When writing the quotients $\varsigma$ in the cochains in the form
\[\varsigma=\sum_{r\in\N}\sum_{k=0}^r\sum_{\E=0}^k\za_{r,k,\E}\frac{x^{\E,k-\E,r-k}}{D},\]
we see that cochains are graded not only with respect to the
cochain-degree $d\in\{0,1,2,3\}$, but also with respect to the
total degree $r$ in $x$ and the partial degree $k$. The preceding
results regarding $X_i(\frac{x^J}{D})$ allow to see that the
coboundary operator is compatible with both degrees, $k$ and $r$,
so that the cohomology can be computed part by part.

\section{Fundamental operators}

Denote by $P_{kr}$ the space of homogeneous polynomials of partial
degree $k$ and total degree $r$, and by $Q_{kr}$ the space
$\frac{1}{D}P_{kr}$. We now study the fundamental operators $X_i$
as endomorphisms of $Q_{kr}$.

It is easy to verify that
\begin{equation}\label{X_1}\begin{array}{ll} &X_1\lp\sum_{\E=0}^k\za_{\E}\frac{x^{\E,k-\E,r-k}}{D}\rp=\\
& \sum_{\E=0}^k [
b(k-\E+1)\za_{\E-1}-a(r-k-1)\za_{\E}-b(\E+1)\za_{\E+1}]\frac{x^{\E,k-\E,r-k}}{D},\end{array}\end{equation}
\begin{equation}\label{X_2}X_2=(2r-3k)\;b\;id_{Q_{kr}},\end{equation} and
\begin{equation}\label{X_3}\begin{array}{ll} &X_3\lp\sum_{\E=0}^k\za_{\E}\frac{x^{\E,k-\E,r-k}}{D}\rp=\\
& \sum_{\E=0}^k [
-2b(k-\E+1)\za_{\E-1}+a(k-2)\za_{\E}+2b(\E+1)\za_{\E+1}]\frac{x^{\E,k-\E,r-k}}{D}.\end{array}\end{equation}

Note also that if $D'=x_1^2+x_2^2$, we have $X_1(D')=0,$
$X_1(x_3)=-ax_3,$ $X_1(D)=-aD,$ $X_2(D')=-2bD',$ $X_2(x_3)=2bx_3,$
$X_2(D)=0,$ $X_3(D')=2aD'$, $X_3(x_3)=0$, $X_3(D)=2aD.$ Hence
$X_1(D^{\E})=-a\E D^{\E},$ $X_3(D^{\E})=2a\E D^{\E},$ for all
$\E\in\Z.$ In particular, $D^{\frac{k}{2}-1}\in
Q_{k,\frac{3}{2}k}$ is an eigenvector of $X_1$ and $X_3$
associated to the eigenvalue $\frac{1}{2}a(2-k)$ and $a(k-2)$
respectively, for all $k\in 2\N.$\\

{\bf Remark} If $\mathbf{Y}$ is $\alpha$, $\alpha Y_1+\beta
Y_2+\gamma Y_3$, $\alpha Y_{23} +\beta Y_{31}+\gamma Y_{12}$, or
$\za Y_{123}$ ($\za,\zb,\zg\in\R$), the cochains
$D^{\E}\mathbf{Y}$ are cocycles for all $\E$ if $a=0$ and for
$\E=0$ otherwise. Indeed, the coboundary of these cochains
vanishes if $[\zL,D^{\E}]=\nabla
D^{\E}$ does.\\

In order to compute the spectrum of the endomorphisms $X_1$ and
$X_3$ of $Q_{kr}$, note that their matrix in the canonical basis
of $Q_{kr}$ is

\[ M_0=\lp\begin{array}{ccccccc}A & B & 0 & \ldots & 0 & 0 & 0\\
-kB & A & 2B & \ldots & 0 & 0 & 0 \\
0 & -(k-1)B & A & \ldots & 0 & 0 & 0 \\
\vdots & & & & & & \vdots\\
0 & 0 & 0 & \ldots & -2B & A & kB \\
0 & 0 & 0 & \ldots & 0  & -B & A
\end{array}\rp,\]
where $(A,B)=(a(k-r+1),-b)$ and $(A,B)=(a(k-2),2b)$ respectively.
A straightforward induction shows that for odd $k$ the determinant
of $M_0-\zl I$ is
\[((A-\zl)^2+B^2)((A-\zl)^2+(3B)^2)\ldots ((A-\zl)^2+(kB)^2),\]
whereas for even $k$ its value is
\[(A-\zl)((A-\zl)^2+(2B)^2)((A-\zl)^2+(4B)^2)\ldots ((A-\zl)^2+(kB)^2).\]
We thus have the
\begin{prop}
(i) For any $k\in 2\N+1$, the operator $X_1$ (resp. $X_3$) has no
eigenvalue.

\noindent (ii) For any $k\in 2\N$, the unique eigenvalue of $X_1$
(resp. $X_3$) is \[\zl = a(k-r+1)\;\;(\mbox{resp.  }a(k-2)).\] The
vector \[\frac{1}{D}(x_1^2+x_2^2)^{\frac{k}{2}}x_3^{r-k}\in
Q_{kr}\] is a basis of eigenvectors.
\end{prop}
Note that this result is an extension of the above remark
regarding eigenvectors of $X_1$ and $X_3$ in the space
$Q_{k,\frac{3}{2}k}$, i.e. for $2r-3k=0.$

\begin{cor}\label{X13}
(i) If $k\in 2\N +1$ and if $k\in 2\N,a\neq 0,k\neq r-1$ (resp.
$k\in 2\N, a\neq 0, k\neq 2$), the operator $X_1$ (resp. $X_3$) is
invertible.

\noindent (ii) If $k\in 2\N\mbox{ and }a=0 \mbox{ or }k=r-1$
(resp. $k\in 2\N\mbox{ and }a=0 \mbox{ or }k=2$), the operator
$X_1$ (resp. $X_3$) is degenerated and the eigenvector
\[\frac{1}{D}(x_1^2+x_2^2)^{\frac{k}{2}}x_3^{r-k}\in Q_{kr}\] is a
basis of the kernel ker$X_1$ (resp. ker$X_3$). Moreover,
\[\mbox{ker}X_1\oplus\mbox{im}X_1=Q_{kr}\mbox{ (resp. }\mbox{ker}X_3\oplus\mbox{im}X_3=Q_{kr}).\]
\end{cor}
{\it Proof.} Only the last result requires an explanation. It
suffices to show that ker$X_i$ $\cap$ im$X_i=0$ ($i\in\{1,3\}$).
Let $X_i(q)\in$ ker$X_i$ ($q\in Q_{kr}$), i.e. $q\in$ ker$X_i^2$.
Since ker$X_i\subset$ ker$X_i^2$, dim ker$X_i^2\ge 1$. We prove by
induction that the eigenvalues of $X_i^2$ are $-((i+1)\E b)^2$
($\E\in\{0,\ldots,\frac{k}{2}\}$), all of them having multiplicity
$2$, except $0$ that has multiplicity $1$. So dim ker$X_i^2=1$ and
ker$X_i=$ ker$X_i^2\ni q.$ \rule{1.5mm}{2.5mm}\\

The following proposition is obvious.
\begin{prop}\label{X2}
(i) If $2r-3k\neq 0$, operator $X_2$ is invertible.

\noindent (ii) If $2r-3k=0,$ we have $X_2=0$ and $X_3=-2X_1.$
\end{prop}

The next proposition is an immediate consequence of the
commutativity of the $Y_i$ ($i\in\{1,2,3\}$).

\begin{prop}\label{co}
(i) All commutators of $X_i$-operators vanish:
\[[X_1,X_2]=[X_1,X_3]=[X_2,X_3]=0.\]

\noindent (ii) If $X_i$ is invertible, \[[X_i^{-1},X_j]=0,\forall
i,j\in\{1,2,3\}.\]
\end{prop}

\section{Injection and short exact sequence}\label{seshds}

Any $1$-cochain $C^1_{kr}$ of degrees $k$ and $r$ can be written
in the form $\frac{1}{D}(p_1Y_1$ $+ p_2Y_2 + p_3Y_3),$ with
$p_i\in P_{kr}$. Conversely, any element of this type reads
$\frac{1}{D}[(p_1x_1-p_2x_2)\p_1+(p_1x_2+p_2x_1)\p_2+p_3x_3\p_3],$
and is induced by a cochain if and only if $p_1x_1-p_2x_2$ is
divisible by $D$ and $p_3$ by $D'=x_1^2+x_2^2$. Indeed, the first
condition implies that $p_1x_2+p_2x_1$ is also a multiple of $D$.
This cochain is then an element of
$P_{k-1,r-2}\p_1+P_{k-1,r-2}\p_2+P_{k-2,r-2}\p_3,$ of course
provided that $k\ge 2,r\ge 2, k\le r-1.$ If $p_3=0$, this
condition reduces to $k\ge 1,r\ge 2, k\le r-1,$ and if
$p_1=p_2=0$, it is replaced by $k\ge 2,r\ge 2, k\le r.$ Of course,
the space \[{\cal P}^1_{kr}=\zD_{kr}Q_{kr}Y_1+ \zD_{kr}Q_{kr}Y_2 +
\zD_{k1}Q_{kr}Y_3\] of {\it potential $1$-cochains} of degrees
$k,r$ (where $\zD_{ij}=1-\zd_{ij},$ defined by means of
Kronecker's symbol, is used in order to group the mentioned cases)
has to be taken into account only if the injected space of {\it
real $1$-cochains} of degrees $k,r$,
\[{\cal R}^1_{kr}=\zD_{kr}P_{k-1,r-2}\p_1+\zD_{kr}P_{k-1,r-2}\p_2+
\zD_{k1}P_{k-2,r-2}\p_3,\] is not vanishing, i.e. if $k\ge 1,r\ge
2, k\le r.$ In the following, we often write $x,y,z$ instead of
$x_1,x_2,x_3$. It is easily checked that the space
\[\begin{array}{l}{\cal S}^1_{kr}=\{\frac{x^{k-1}z^{r-k}}{D}[\zD_{kr}exY_1\\+
\zD_{kr}fxY_2+\zD_{k1}(gx+hy)Y_3]: c,d,e,f\in\R\}\end{array}\] is
supplementary to ${\cal R}^1_{kr}$ in ${\cal P}^1_{kr}.$ Similar
spaces ${\cal P}^d_{kr},{\cal R}^d_{kr},{\cal S}^d_{kr}$ can be
defined for $d=0,k\ge 2,r\ge 3,k\le r-1;$ $d=2,k\ge 0,r\ge 1,k\le
r$ and $d=3,k\ge 0,r\ge0,k\le r.$ These spaces are described at
the end of this section. Hence the whole space of potential
cochains ${\cal P}=\oplus_{d,k,r}{\cal P}^d_{kr}$ is the direct
sum of the whole space of real cochains ${\cal
R}=\oplus_{d,k,r}{\cal R}^d_{kr}$ and the supplementary space
${\cal S}=\oplus_{d,k,r}{\cal S}^d_{kr}$ (we can view $k$ and $r$
as subscripts running through $\N$ ($k\le r$), provided that the
spaces associated to forbidden values are considered as
vanishing).\\

The spaces $({\cal P},\p_{\cal P})$ and $({\cal R},\p_{\cal R})$,
where $\p_{\cal P}=\p_{\cal R}=[\zL,\cdot]$, are differential
spaces. Denote by $p_{\cal R}$ and $p_{\cal S}$ the projections of
${\cal P}$ onto ${\cal R}$ and ${\cal S}$ respectively and set for
any $s\in{\cal S}$,
\[\zf s=p_{\cal R}\p_{\cal P}s,\p_{\cal S}s=p_{\cal S}\p_{\cal
P}s.\]

\begin{prop}
(i) The endomorphism $\p_{\cal S}\in End{\cal S}$ is a
differential on ${\cal S}$.

\noindent (ii) The linear map $\zf\in {\cal L}({\cal S},{\cal R})$
is an anti-homomorphism of differential spaces from $({\cal
S},\p_{\cal S})$ into $({\cal R},\p_{\cal R})$.
\end{prop}
{\it Proof.} Direct consequence of $\p_{\cal P}^2=0.$
\rule{1.5mm}{2.5mm}

\begin{prop}
If $i$ denotes the injection of ${\cal R}$ into ${\cal P}$, the
sequence
\[0\raa {\cal R}\stackrel{i}{\raa}{\cal P}\stackrel{p_{\cal S}}{\raa}{\cal S}\raa 0\]
is a short exact sequence of homomorphisms of differential spaces.
Hence the triangle

\begin{center}
\begin{picture}(100,100)(12.5,0)
\put(50,20){\vector(-1,2){30}} \put(5,85){$H({\cal R})$}
\put(35,87.5){\vector(1,0){60}}\put(62.5,97.5){\makebox(0,0){$i_{\sharp}$}}
\put(100,85){$H({\cal P})$}\put(110,80){\vector(-1,-2){30}}
\put(52.5,10){$H({\cal
S})$}\put(20,50){$\zf_{\sharp}$}\put(105,50){$(p_{\cal
S})_{\sharp}$}
\end{picture}
\end{center}
is exact.\end{prop} {\it Proof.} We only need check that the
linear map $\zf_{\sharp}$ induced by $\zf$ coincides with the
connecting homomorphism $\zD$. It suffices to remember the
definition of $\zD$. \rule{1.5mm}{2.5mm}\\

Remark now that the exact triangle induces for any
$(k,r)\in\N^2,k\le r$, a "long" exact sequence of linear maps:
\[\begin{array}{l}0\raa H^0_{kr}({\cal
R})\stackrel{i_{\sharp}}{\raa}\ldots\\\;\;\;\;\;\;\;
\stackrel{\zf_{\sharp}}{\raa} H^d_{kr}({\cal
R})\stackrel{i_{\sharp}}{\raa}H^d_{kr}({\cal P})
\stackrel{(p_{\cal S})_{\sharp}}{\raa}H^d_{kr}({\cal
S})\stackrel{\zf_{\sharp}}{\raa}H^{d+1}_{kr}({\cal
R})\stackrel{i_{\sharp}}{\raa}\ldots\stackrel{(p_{\cal
S})_{\sharp}}{\raa}H^3_{kr}({\cal S})\raa 0.\end{array}\] We
denote the kernel and the image of the restricted map
$\zf_{\sharp}\in {\cal L}(H^d_{kr}({\cal S}),$ $H^{d+1}_{kr}({\cal
R}))$ by ker$^d_{kr}\zf_{\sharp}\subset H^d_{kr}({\cal S})$ and
im$^{d+1}_{kr}\zf_{\sharp}\subset H^{d+1}_{kr}({\cal R}).$ Similar
notations are used if $i_{\sharp}$ or $(p_{\cal S})_{\sharp}$ are
viewed as restricted maps.
\begin{cor}\label{dirsum}
For any $d\in\{0,1,2,3\},k\in\N,r\in\N,k\le r$, we
have\[\begin{array}{ll}H^d_{kr}({\cal R})& \simeq
\mbox{ker}^d_{kr}i_{\sharp}\oplus\mbox{im}^d_{kr}i_{\sharp}\\& =
\mbox{im}^d_{kr}\zf_{\sharp}\oplus\mbox{ker}^d_{kr}(p_{\cal
S})_{\sharp}\\& \simeq H^{d-1}_{kr}({\cal
S})/ker^{d-1}_{kr}\zf_{\sharp}\oplus H^d_{kr}({\cal P
})/ker^d_{kr}\zf_{\sharp}.\end{array}\]
\end{cor}
{\it Proof.} Apply exactness of the long sequence.
\rule{1.5mm}{2.5mm}\\

We will compute the ${\cal R}$-cohomology by computing the simpler
${\cal S}$- and ${\cal P}$-cohomology (and in some cases the
anti-homomorphism $\zf$).

The preceding result is easily understood. The space $Z_{\cal R}$
of ${\cal R}$-cocycles is a subset of $Z_{\cal P}$. Among the
${\cal P}$-classes there may be classes without representatives in
$Z_{\cal R}$. Take now the classes with a non-empty intersection
with $Z_{\cal R}$. Two cocycles in different intersections can not
be ${\cal R}$-cohomologous. Two cocycles in the same intersection
are or not ${\cal R}$-cohomologous. Hence the picture of the
${\cal P}$- and ${\cal R}$-cohomologies. Remember now that the
isomorphism $H({\cal R})\simeq ker i_{\sharp}\oplus im i_{\sharp}$
means that $H({\cal R})= ker i_{\sharp}\oplus \chi(im
i_{\sharp})$, where $\chi$ is a linear right inverse of
$i_{\sharp}$. The space $im i_{\sharp}=\{[\zr]_{\cal P}:\zr\in
Z_{\cal R}\}$, where $[\cdot]_{\cal P}$ denotes a class in ${\cal
P}$, is the space of ${\cal P}$-classes with intersection. The
space $\chi(im i_{\sharp})\simeq im i_{\sharp}$ is made up of one
of the source classes of each ${\cal P}$-class with intersection.
We obtain the missing ${\cal R}$-classes when adding the kernel.

As for the meaning of $\zf$, since $\p_{\cal S}s=\p_{\cal P}s-\zf
s,$ this anti-homomorphism is nothing but the correction that
turns the Poisson-coboundary in a coboundary on ${\cal S}$. If
$\zf$ were $0$, we would of course have $H({\cal R})=H({\cal
P})/H({\cal S})$. We recover this result as special case of the
preceding corollary.

\begin{cor}\label{HS0}
For any $d\in\{0,1,2,3\},k\in\N,r\in\N,k\le r$,\\

\noindent (i) if
\[H^{d-1}_{kr}({\cal S})=0 \;\;\;(\mbox{resp. } H^{d-1}_{kr}({\cal
S})=H^d_{kr}({\cal S})=0),\] then \[i_{\sh}\in Isom(H^d_{kr}({\cal
R}),im^d_{kr}i_{\sh}) \;\;\;(\mbox{resp. } i_{\sh}\in
Isom(H^d_{kr}({\cal R}),H^d_{kr}({\cal P}))),\]

\noindent (ii) if \[H^{d}_{kr}({\cal P})=0 \;\;\;(\mbox{resp. }
H^{d-1}_{kr}({\cal P})=H^d_{kr}({\cal P})=0),\] then
\[\zf_{\sh}\in Isom(H^{d-1}_{kr}({\cal S})/ker^{d-1}_{kr}\zf_{\sh},H^d_{kr}({\cal R}))
\;\;\;(\mbox{resp. } \zf_{\sh}\in Isom(H^{d-1}_{kr}({\cal
S}),H^d_{kr}({\cal R}))).\]
\end{cor}

\noindent {\it Proof.} Obvious.  \rule{1.5mm}{2.5mm}\\


The below basic formulas are obtained by straightforward
computations. For instance, a potential $0$-cochain of degree
$(k,r)$, $\pi=\frac{p}{D}\in{\cal P}^0_{kr}=Q_{kr},$ is a member
of the corresponding real cochain space ${\cal
R}^0_{kr}=P_{k-2,r-3}$ ($k\ge 2,r\ge 3,k\le r-1$), if and only if
$p\in P_{kr}$ is divisible by $D'=x_1^2+x_2^2=x^2+y^2$. If
$p=\sum_{\E=0}^k \za_{\E} x^{\E,k-\E,r-k}$, this divisibility
condition means that $\za_0-\za_2+\za_4-\ldots=0$ and
$\za_1-\za_3+\za_5-\ldots=0$. Hence, a potential cochain can be
made a real cochain by changing the coefficients $\za_{k-1}$ and
$\za_k$, so that any potential cochain can be written in a unique
way as the sum of a real cochain and an element of ${\cal
S}^0_{kr}=\{\frac{x^{k-1}z^{r-k}}{D}(cx+dy):c,d\in\R\}$.

\begin{center}{\bf Formulary 1}\end{center}

\noindent {\bf 1. 0-cochains}\\

\noindent $k\ge 2,r\ge 3,k\le r-1$:\\

\noindent ${\cal P}^0_{kr}=Q_{kr}$\\
\noindent ${\cal R}^0_{kr}=P_{k-2,r-3}$\\
\noindent ${\cal S}^0_{kr}=\{\frac{x^{k-1}z^{r-k}}{D}(cx+dy):c,d\in\R\}$\\

\noindent A potential cochain $\pi=\frac{p}{D}$ is a real cochain
if and only if $p$ is divisible by $D'$\\

\noindent {\bf 2. 1-cochains}\\

\noindent $k\ge 1,r\ge 2,k\le r$:\\

\noindent ${\cal P}^1_{kr}=\zD_{kr}Q_{kr}Y_1+ \zD_{kr}Q_{kr}Y_2 +
\zD_{k1}Q_{kr}Y_3$\\

\noindent ${\cal
R}^1_{kr}=\zD_{kr}P_{k-1,r-2}\p_1+\zD_{kr}P_{k-1,r-2}\p_2+
\zD_{k1}P_{k-2,r-2}\p_3$\\

\noindent ${\cal
S}^1_{kr}=\{\frac{x^{k-1}z^{r-k}}{D}[\zD_{kr}exY_1+
\zD_{kr}fxY_2+\zD_{k1}(gx+hy)Y_3]:\\e,f,g,h\in\R\}$\\

\noindent A potential cochain
$\pi=\frac{1}{D}[\zD_{kr}p_1Y_1+\zD_{kr}p_2Y_2+\zD_{k1}p_3Y_3]$ is
a real cochain if and only if $\zD_{kr}[p_1x_1-p_2x_2]$ and
$\zD_{k1}p_3$ are divisible by $D'$\\

\noindent {\bf 3. 2-cochains}\\

\noindent $k\ge 0,r\ge 1,k\le r$:\\

\noindent ${\cal
P}^2_{kr}=\zD_{k0}Q_{kr}Y_{23}+\zD_{k0}Q_{kr}Y_{31}+\zD_{kr}Q_{kr}Y_{12}$\\

\noindent ${\cal
R}^2_{kr}=\zD_{k0}P_{k-1,r-1}\p_{23}+\zD_{k0}P_{k-1,r-1}\p_{31}+\zD_{kr}P_{k,r-1}\p_{12}$\\

\noindent ${\cal
S}^2_{kr}=\{\frac{x^kz^{r-k}}{D}[\zD_{k0}iY_{23}+\zD_{k0}jY_{31}]:i,j\in\R\}$\\

\noindent A potential cochain
$\pi=\frac{1}{D}[\zD_{k0}p_1Y_{23}+\zD_{k0}p_2Y_{31}+\zD_{kr}p_3Y_{12}]$
is a real cochain if and only if $\zD_{k0}[p_1x_1-p_2x_2]$ is
divisible by $D'$\\

\noindent {\bf 4. 3-cochains}\\

\noindent $k\ge 0,r\ge 0,k\le r$:\\

\noindent ${\cal P}^3_{kr}=Q_{kr}Y_{123}$\\

\noindent ${\cal R}^3_{kr}=P_{kr}\p_{123}$\\

\noindent ${\cal S}^3_{kr}=0$\\

\noindent The spaces of potential and real cochains coincide\\

We write the coefficients $(e,f,g,h)$ (resp. $(i,j)$) of the
coboundary $\p_{\cal S}\zs^0_{kr}$ (resp. $\p_{\cal S}\zs^1_{kr})$
in terms of the coefficients $(c,d)$ (resp.
$(e\zD_{kr},f\zD_{kr},g\zD_{k1},h\zD_{k1})$) of the supplementary
cochain $\zs^0_{kr}\in {\cal S}^0_{kr},$ $k\ge 2,r\ge 3,k\le r-1$
(resp. $\zs^1_{kr}\in {\cal S}^1_{kr},$ $k\ge 1,r\ge 2,k\le r$)
and of the Pauli type matrices
\[\zm_0=\lp\begin{array}{cc}1&0\\0&1\end{array}\rp,\zm_1=\lp\begin{array}{cc}0&1\\1&0\end{array}\rp,
\zm_2=\lp\begin{array}{cc}0&1\\-1&0\end{array}\rp,\zm_3=\lp\begin{array}{cc}1&0\\0&-1\end{array}\rp.\]
The expressions of $\zf(\zs^d_{kr})$ are not all indispensable.

\begin{center}{\bf Formulary 2}\end{center}

\noindent {\bf 1. 0-cochains}\\

\noindent Coefficients of $\p_{\cal S}\zs^0_{kr}$ in terms of the
coefficients of $\zs^0_{kr},$ $k\ge 2,r\ge 3,k\le r-1$:
\[\lp\begin{array}{c}e\\f\\\hline g\\h\end{array}\rp=\lp\begin{array}{c}2b(r-k)\zm_1-a(r-k-1)\zm_3\\\hline
a(k-2)\zm_0-2bk\zm_2\end{array}\rp\lp\begin{array}{c}c\\d\end{array}\rp\]

\noindent {\bf 2. 1-cochains}\\

\noindent Coefficients of $\p_{\cal S}\zs^1_{kr}$ in terms of the
coefficients of $\zs^1_{kr},$ $k\ge 1,r\ge 2,k\le r$:
\[\lp\begin{array}{c}i\\j\end{array}\rp=\lp\begin{array}{ccccc}2bk\zm_0-a(k-2)\zm_2
\mid a(r-k-1)\zm_1+2b(r-k)\zm_3\end{array}\rp
\lp\begin{array}{c}e\zD_{kr}\\f\zD_{kr}\\\hline
g\zD_{k1}\\h\zD_{k1}\end{array}\rp\] \noindent Value of
$\zf(\zs^1_{rr})$ in terms of the coefficients of $\zs^1_{rr},$
$r\ge 2,a=0$:
\[\begin{array}{lll}\zf(\zs^1_{rr})&=&\frac{bx^{r-2}}{D}[-rx(gx+hy)Y_{23}+(-h\,x^2+rg\,xy+(r-1)h\,y^2)
Y_{31}]\\&=&-bx^{r-2}[(rgx+(r-1)hy)\p_{23}+hx\p_{31}]\end{array}\]

\noindent {\bf 3. 2-cochains}\\

\noindent Coboundary $\p_{\cal S}\zs^2_{kr},$ $k\ge 0,r\ge 1,k\le
r$:
\[\p_{\cal S}\zs^2_{kr}=0\]
\noindent Value of $\zf(\zs^2_{rr})$ in terms of the coefficients
of $\zs^2_{rr},$ $r\ge 1$:
\begin{equation}\label{f}\zf(\zs^2_{rr})=\frac{x^{r-1}}{D}[(ai-brj)x-briy]Y_{123}\end{equation}

\noindent {\bf 4. 3-cochains}\\

\noindent Coboundary $\p_{\cal S}\zs^3_{kr},$ $k\ge 0,r\ge 0,k\le
r$:
\[\p_{\cal S}\zs^3_{kr}=0\]

\section{0 - cohomology spaces}

\subsection{$\mathbf{{\cal S}}$ - cohomology}

\begin{prop}
The $0$ - cohomology space of ${\cal S}$ vanishes: $H^0({\cal
S})=0.$
\end{prop}

\noindent {\it Proof.} If $(c,d)$ are the coefficients of an
arbitrary cochain $\zs^0_{kr}$ ($k\ge 2,r\ge 3,k\le r-1$), we have
for instance
\[M_0\lp\begin{array}{c}c\\d\end{array}\rp=0,\;\;\;M_0=a(k-2)\zm_0-2bk\zm_2.\]
Since $det\,M_0=a^2(k-2)^2+4b^2k^2>0,$ it follows that
$\zs^0_{kr}=0.$  \rule{1.5mm}{2.5mm}\\

\begin{cor}\label{dimSim1} For any $k\ge 2,r\ge 3,k\le r-1$, we have
$dim\;im^1_{kr}\;\p_{\cal S}=dim\;{\cal S}^0_{kr}=2.$\end{cor}

\subsection{$\mathbf{{\cal P}}$ - cohomology and $\mathbf{{\cal R}}$ - cohomology}

\begin{theo}
The $0$ - cohomology spaces of ${\cal P}$ and ${\cal R}$
coincide:\\
(i) if $a\neq 0$, $H^0({\cal P})=H^0({\cal R})=H^0_{23}({\cal P})=H^0_{23}({\cal R})=\R,$\\
(ii) if $a=0$, $H^0({\cal P})=H^0({\cal
R})=\bigoplus_{k=1}^{\infty}H^0_{2k,3k}({\cal
P})=\bigoplus_{k=1}^{\infty}H^0_{2k,3k}({\cal
R})=\bigoplus_{m=0}^{\infty}\R\, D^m.$
\end{theo}

\noindent {\it Proof.} The equality of both cohomologies is a
consequence of Corollary \ref{HS0}. So let $k\ge 2,r\ge 3,k\le
r-1$ and $\zp^0_{kr}=q\in {\cal P}^0_{kr}\cap ker\p_{\cal P}$. We
have $X_1(q)=X_2(q)=0$. Apply now Proposition \ref{X2} and
Corollary \ref{X13}. If $3k-2r\neq 0$ and if $3k-2r=0,$ $a\neq 0$
and $k\neq r-1$, the cocycle vanishes. If $3k-2r=0$ and $a=0$ or
$k=r-1$, $q=\za D^{\frac{k}{2}-1}$ ($\za\in\R$).
\rule{1.5mm}{2.5mm}

\section{1 - cohomology spaces}

\subsection{$\mathbf{{\cal S}}$ - cohomology}

\begin{prop}\label{H1S} If $a\neq 0$, $H^1({\cal S})=0,$ and if $a=0$, $H^1({\cal S})=
\bigoplus_{m=2}^{\infty}H^1_{mm}({\cal
S})=\bigoplus_{m=2}^{\infty}\frac{x_1^{m-1}}{D}(\R x_1+\R
x_2)Y_3.$
\end{prop}

\noindent {\it Proof.} Set $M_1=2bk\zm_0-a(k-2)\zm_2$ and
$N_1=a(r-k-1)\zm_1+2b(r-k)\zm_3$ ($k\ge 1,r\ge 2,k\le r$). The
coefficients of any $1$-cocycle of degrees $k,r$ verify
\[\zD_{kr}\lp\begin{array}{c} e\\
f\end{array}\rp=-\zD_{k1}M_1^{-1}N_1\lp\begin{array}{c}g\\h\end{array}\rp.\]
If $k=1$, the cocycle vanishes. If $k=r$ the cocycle equation
reads $N_1\lp\begin{array}{c} g\\
h\end{array}\rp=0$, with $N_1=-a\zm_1$. For $a\neq 0$, the cocycle
vanishes again, and for $a=0$, the cohomology space is
$H^1_{rr}({\cal S})=\frac{x^{r-1}}{D}(\R x+\R y)Y_3$, since
coboundaries do not exist for $k=r$. Finally, if $2\le k\le r-1,
r\ge 3,$ Corollary \ref{dimSim1} shows that there is no
cohomology. (A direct proof of this result, based upon the
multiplication table
$\zm_i\zm_j=(-1)^{ij+1}[\zd_{ij}\zm_0+\ze_{ijk}\zm_k]$
($i,j\in\{1,2,3\}$), where $\ze_{ijk}$ is the Levi-Civita symbol,
is also easy.) \rule{1.5mm}{2.5mm}


\subsection{$\mathbf{{\cal P}}$ - cohomology and $\mathbf{{\cal R}}$ - cohomology}

\begin{theo}\label{H1PR} The first cohomology groups of ${\cal P}$ and ${\cal R}$ are isomorphic:\\
(i) $H^1({\cal P})=H^1({\cal R})=H^1_{23}({\cal P})=H^1_{23}({\cal
R})
=\R Y_1+\R Y_2+\R Y_3$, if $a\neq 0$,\\
(ii) $H^1({\cal P})=H^1({\cal
R})=\bigoplus_{k=1}^{\infty}H^1_{2k,3r}({\cal
P})=\bigoplus_{k=1}^{\infty}H^1_{2k,3r}({\cal
R})=\bigoplus_{m=0}^{\infty}D^{m}(\R Y_1+\R Y_2+\R Y_3),$ if
$a=0.$
\end{theo}

\noindent {\it Proof.} Let
$\zp^1_{kr}=\zD_{kr}q_1Y_1+\zD_{kr}q_2Y_2+\zD_{k1}q_3Y_3$ ($k\ge
1,r\ge 2,k\le r$) be a member of ${\cal P}^1_{kr}\cap ker\p_{\cal
P}$. The cocycle equation reads
\[\begin{array}{l}X_2(\zD_{k1}q_3)-X_3(\zD_{kr}q_2)=0,\\
X_3(\zD_{kr}q_1)-X_1(\zD_{k1}q_3)=0,\\
X_1(\zD_{kr}q_2)-X_2(\zD_{kr}q_1)=0.\end{array}\] Moreover, for
$k\ge 2,r\ge 3,k\le r-1$, this cocycle can be a coboundary, i.e.
there is $\zp^0_{kr}=q$ in ${\cal P}^0_{kr}$ such that
\[X_1(q)=q_1,\;X_2(q)=q_2,\;X_3(q)=q_3.\]
{\bf 1.  $\mathbf{k=1}$ or $\mathbf{k=r}$}\\

\noindent We treat the first case (resp. the second case). The
cocycle equation implies $X_3(q_1)=X_3(q_2)=0$ (resp.
$X_2(q_3)=0$). In view of Corollary \ref{X13} (resp. Proposition
\ref{X2}) $X_3$
(resp. $X_2$) is invertible and cocycles vanish.\\

\noindent {\bf 2.  $\mathbf{2\le k\le r-1}$} (then $r\ge 3$)\\

\noindent {\bf 2.a.  $\mathbf{2r-3k\neq 0}$}\\

\noindent Set $q=X_2^{-1}(q_2)$. Due to Proposition \ref{co} and
the cocycle relations, we then have $X_1(q)=q_1,X_3(q)=q_3,$ so
that all cocycles are coboundaries.\\

\noindent {\bf 2.b.  $\mathbf{2r-3k=0}$} (then $k\in 2\N$)\\

\noindent Here $X_2=0$ and $X_3=-2X_1$. The cocycle condition is
$X_1(q_2)=0$ and $X_3(q_1)=X_1(q_3),$ i.e. $X_1(2q_1+q_3)=0$.\\

\noindent {\bf 2.b.I  $\mathbf{a\neq 0}$ and $\mathbf{k\neq r-1}$}
(i.e. $a\neq 0$ and $(k,r)\neq (2,3)$)\\

\noindent In this case $q_2=0$ and we choose $q=X_1^{-1}(q_1)$.
This entails that $X_3(q)=q_3.$\\

\noindent {\bf 2.b.II  $\mathbf{a=0}$ or $\mathbf{k=r-1}$}
(i.e. $a=0$ or $(k,r)=(2,3)$)\\

\noindent Corollary \ref{X13} allows to decompose $q_1$ in the
form $q_1=\za_1D^{\frac{k}{2}-1}+X_1(q)$ ($\za_1\in\R,q\in
Q_{kr}$). It is clear that $q_2=\za_2D^{\frac{k}{2}-1}+X_2(q)$
($\za_2\in\R$) and that
$q_3=\za_3'D^{\frac{k}{2}-1}-2q_1=\za_3D^{\frac{k}{2}-1}+X_3(q)$
($\za_3',\za_3\in\R$). Hence, the considered cocycle is
cohomologous to $D^{\frac{k}{2}-1}(\za_1Y_1+\za_2Y_2+\za_3Y_3)$.
As the kernel and the image of $X_1$ and $X_3$ are supplementary
and $X_2$ vanishes, such cocycles---with different
$\za$-coefficients---can not be cohomologous.\\

As for the isomorphism between the ${\cal P}$- and ${\cal
R}$-cohomology, remember that $H^0({\cal S})=0$. Corollary
\ref{HS0} then implies that $i_{\sh}\in Isom(H^1_{kr}({\cal
R}),H^1_{kr}({\cal P}))$, for all $k,r$. Indeed, if
$H^1_{kr}({\cal S})\neq0$, we have $k=r\ge 2$ and $H^1_{kr}({\cal
P})=0=im^1_{kr}i_{\sh}$.  \rule{1.5mm}{2.5mm}

\section{2 - cohomology spaces}

\subsection{$\mathbf{{\cal S}}$ - cohomology}

\begin{prop} If $a\neq 0,$ $H^2({\cal S})=H^2_{11}({\cal S})=\frac{x_1}{D}(\R Y_{23}+\R
Y_{31}),$ and if $a=0,$ $H^2({\cal
S})=\bigoplus_{m=1}^{\infty}H^2_{mm}({\cal
S})=\bigoplus_{m=1}^{\infty}\frac{x_1^m}{D}(\R Y_{23}+\R Y_{31}).$
\end{prop}

\noindent {\it Proof.}  As ${\cal S}^3=0,$ any cochain is a
cocycle. Let $\zs^2_{kr}$ ($k\ge 0,r\ge 1,k\le r$) be a cocycle
with coefficients $(\zD_{k0}i,\zD_{k0}j)$. If $k\ge 1,r\ge 2,k\le
r$, it is a coboundary if there are coefficients
$(\zD_{kr}e,\zD_{kr}f,\zD_{k1}g,\zD_{k1}h)$, such that
\[\lp\begin{array}{c}i\\j\end{array}\rp=\lp\begin{array}{ccccc}2bk\zm_0-a(k-2)\zm_2\mid
a(r-k-1)\zm_1+2b(r-k)\zm_3\end{array}\rp
\lp\begin{array}{c}\zD_{kr}e\\\zD_{kr}f\\\hline
\zD_{k1}g\\\zD_{k1}h\end{array}\rp.\]

\noindent If $r\ge 1, k=0$, all cochains vanish. If $r=1, k=1$,
there are no coboundaries. Consider now the case $r\ge 2, k\ge 1$.
If $k\neq r$ and if $k=r$ and $a\neq 0$, any cocycle is a
coboundary. If $k=r$ and $a=0$, the unique coboundary is $0$.
\rule{1.5mm}{2.5mm}

\subsection{$\mathbf{{\cal P}}$ - cohomology}

\begin{prop}\label{H2P}
The second cohomology group of the complex $({\cal P},\p_{\cal
P})$ is,\\(i) if $a\neq 0$,\\\[\begin{array}{lll}H^2({\cal P})&=&H^2_{23}({\cal
P})\oplus H^2_{11}({\cal P})
\\&=&\R
Y_{23}+\R Y_{31}+\R Y_{12}\\&&\oplus\frac{1}{D}\left[\R\lp
(D-\frac{a}{b}x_1x_2x_3)\p_{23}+\frac{a}{b}x_1^2x_3\p_{31}\rp
+\R\lp(D+\frac{a}{b}x_1x_2x_3)\p_{31}-\frac{a}{b}x_2^2x_3\p_{23}\rp\right],\end{array}\]
(ii) if $a=0$, \\\[\begin{array}{lll}H^2({\cal P})&=&\bigoplus_{k=1}^{\infty}H^2_{2k,3r}({\cal P})\oplus H^2_{11}({\cal
P})\oplus\bigoplus_{m=1}^{\infty}H^2_{0m}({\cal P})
\\&=&\bigoplus_{m=0}^{\infty}D^{m}\lp\R
Y_{23}+\R Y_{31}+\R Y_{12}\rp\\&&\oplus\;
\R\p_{23}+\R\p_{31}\oplus\bigoplus_{m=0}^{\infty}\R
x_3^{m}\p_{12}.\end{array}\]
\end{prop}
\noindent {\it Proof.}  The cocycle equation for a cochain
$\zp^2_{kr}=\zD_{k0}q_1Y_{23}+\zD_{k0}q_2Y_{31}+\zD_{kr}q_3Y_{12}$
($k\ge 0,r\ge 1,k\le r;q_1,q_2,q_3\in Q_{kr}$) reads
\[X_1(\zD_{k0}q_1)+X_2(\zD_{k0}q_2)+X_3(\zD_{kr}q_3)=0.\]
This cocycle can be the coboundary of a cochain
$\zp^1_{kr}=\zD_{kr}\zvr_1Y_1+\zD_{kr}\zvr_2Y_2+\zD_{k1}\zvr_3Y_3$
($k\ge 1,r\ge 2,k\le r;\zvr_1,\zvr_2,\zvr_3\in Q_{kr}$):
\[\begin{array}{l}X_2(\zD_{k1}\zvr_3)-X_3(\zD_{kr}\zvr_2)=\zD_{k0}q_1,\\
X_3(\zD_{kr}\zvr_1)-X_1(\zD_{k1}\zvr_{3})=\zD_{k0}q_2,\\
X_1(\zD_{kr}\zvr_2)-X_2(\zD_{kr}\zvr_1)=\zD_{kr}q_3.\end{array}\]
{\bf 1. $\mathbf{k=0,r\ge 1}$} (then $k\le r$)\\

\noindent We get $X_3(q_3)=0.$ If $a\neq 0,$ the cocycle
$\zp^2_{0r}$ vanishes, and if $a=0,$ it has the form
\[\zp^2_{0r}=\frac{\za_0}{D}z^rY_{12}=\za_0z^{r-1}\p_{12}\;\;\;(a=0,r\ge 1,\za_0\in\R).\]

\noindent {\bf 2. $\mathbf{k\ge 1,r\ge 1,k=r}$}\\

\noindent Cocycle condition: $q_2=-X_2^{-1}X_1(q_1)$.\\

\noindent {\bf 2.a. $\mathbf{k=r=1}$}\\

\noindent It follows immediately from Equation \ref{X_1}, Equation
\ref{X_2} and Equation \ref{Y}, that
\begin{equation}\label{pi211}\begin{array}{lll}\zp^2_{11}&=&\frac{1}{D}(\za_1x+\za_2y)Y_{23}+\frac{1}{D}
\lp(\za_2+\frac{a}{b}\za_1)x+(\frac{a}{b}\za_2-\za_1)y\rp
Y_{31}\\&=&\frac{1}{D}\left[\za_1\lp(D-\frac{a}{b}xyz)\p_{23}+\frac{a}{b}x^2z\p_{31}\rp
+\za_2\lp
(D+\frac{a}{b}xyz)\p_{31}-\frac{a}{b}y^2z\p_{23}\rp\right]
\;\;\;(\za_1,\za_2\in\R).\end{array}\end{equation}

\noindent {\bf 2.b. $\mathbf{k=r\ge 2}$}\\

\noindent As the coboundary condition reads
$X_2(\zvr_3)=q_1,X_1(\zvr_3)=-q_2,$ we only need choose
$\zvr_3=X_2^{-1}(q_1).$\\

\noindent {\bf 3. $\mathbf{k\ge 1,k\le r-1}$} (then $r\ge 2$)\\

\noindent Cocycle condition: $X_1(q_1)+X_2(q_2)+X_3(q_3)=0$.\\

\noindent {\bf 3.a. $\mathbf{k=1}$}\\

\noindent Coboundary condition:
$X_3(\zvr_2)=-q_1,X_3(\zvr_1)=q_2,X_1(\zvr_2)-X_2(\zvr_1)=q_3.$ It
suffices to take $\zvr_1=X_3^{-1}(q_2),\zvr_2=-X_3^{-1}(q_1).$\\

\noindent {\bf 3.b. $\mathbf{k\ge 2}$} (then $r\ge 3$)\\

\noindent All $\zD$-factors in the cocycle and coboundary
conditions disappear. \\

\noindent {\bf 3.b.I. $\mathbf{3k-2r\neq 0}$}\\

\noindent Set
$\zvr_2=0,\zvr_3=X_2^{-1}(q_1),\zvr_1=-X_2^{-1}(q_3).$\\

\noindent {\bf 3.b.II. $\mathbf{3k-2r=0}$}\\

\noindent Since $X_2=0$ and $X_3=-2X_1$, the equations read
$X_1(q_1-2q_3)=0$ and
$X_1(\zvr_2)=\frac{1}{2}q_1,X_1(2\zvr_1+\zvr_3)=-q_2,X_1(\zvr_2)=q_3$
respectively.\\

\noindent {\bf 3.b.II.}$\mathbf{\alpha}$. $\mathbf{a\neq
0,(k,r)\neq
(2,3)}$\\

\noindent As $X_1$ is invertible, $q_1=2q_3$ and we set
$\zvr_1=0,\zvr_2=X_1^{-1}(q_3),\zvr_3=-X_1^{-1}(q_2).$\\

\noindent {\bf 3.b.II.$\mathbf{\beta}$. $\mathbf{a=0}$ or
$\mathbf{(k,r)=
(2,3)}$}\\

\noindent Corollary \ref{X13} allows to write $q_i$
($i\in\{1,2,3\}$) in a unique way as the sum
$q_i=\za_{2+i}D^{\frac{k}{2}-1}+q_i'$, $\za_{2+i}\in\R,q_i'\in
Q_{k,\frac{3}{2}k}$ of an element of $ker\,X_1$ and an element of
$im\,X_1$. As
$q_1-2q_3-(\za_3-2\za_5)D^{\frac{k}{2}-1}=q_1'-2q_3'\in
ker\,X_1\cap im\,X_1$, we conclude that $q_1'=2q_3'$. Hence
$q_1'Y_{23}+q_2'Y_{31}+q_3'Y_{12}$ is the coboundary of
$\zvr_1Y_1+\zvr_2Y_2+\zvr_3Y_3,$ where
$\zvr_1=0,X_1(\zvr_2)=q_3',X_1(\zvr_3)=-q_2'.$ Finally, the
original cocycle is cohomologous to
\[\begin{array}{c}\zp^2_{k,\frac{3}{2}k}\simeq
D^{\frac{k}{2}-1}\lp\za_3Y_{23}+\za_{4}Y_{31}+\za_{5}Y_{12}\rp\\\;\;\;(a=0\mbox{
and }k\in\{2,4,6,\ldots\}\mbox{  or  }a\neq 0\mbox{ and }k=2;
\za_3,\za_4,\za_5\in\R).\end{array}\] Due to the supplementary
character of the kernel and the image of $X_1$, two different
cocycles of this type can not be cohomologous. \rule{1.5mm}{2.5mm}

\subsection{$\mathbf{{\cal R}}$ - cohomology}

\begin{theo}(i) If $a\neq 0$,\\\[H^2({\cal R})=H^2_{23}({\cal R})=\R
Y_{23}+\R Y_{31}+\R Y_{12},\]
(ii) if $a=0$, \\\[\begin{array}{lll}H^2({\cal R})&=&\bigoplus_{k=1}^{\infty}H^2_{2k,3k}({\cal
R})\oplus\bigoplus_{m=1}^{\infty}H^2_{mm}({\cal
R})\oplus\bigoplus_{m=1}^{\infty}H^2_{0m}({\cal R})
\\&=&\bigoplus_{m=0}^{\infty}D^{m}\lp\R
Y_{23}+\R Y_{31}+\R Y_{12}\rp\\&&\oplus\;
\bigoplus_{m=1}^{\infty}x_1^{m-2}\lp\R
x_1\p_{23}+\R(x_1\p_{31}+(m-1)x_2\p_{23})\rp\oplus\bigoplus_{m=0}^{\infty}\R
x_3^{m}\p_{12}.\end{array}\]
\end{theo}

\noindent {\it Proof.}  We apply once more Corollary \ref{HS0}.

If $(kr)\neq (mm)$ ($m\ge 1$), the map $i_{\sh}$ is an isomorphism
between $H^2_{kr}({\cal R})$ and $H^2_{kr}({\cal P})$.

If $(kr)=(11),$ we get $i_{\sh}\in Isom(H^2_{11}({\cal
R}),im^2_{11}\,i_{\sh})$. Since no coboundaries exist for this
degree, the ${\cal R}$-cohomology is made up by the ${\cal
P}$-cocycles in ${\cal R}$. But, $\pi^2_{11}\in {\cal R}$, see
\ref{pi211}, if and only if $a\za_1=a\za_2=0$. So, for $a\neq0$,
we have $H^2_{11}({\cal R})=0,$ whereas for $a=0,$ we obtain
$H^2_{11}({\cal R})=\R\p_{23}+\R\p_{31}$.

Consider now $(kr)=(mm)$ ($m\ge 2$). As $H^1_{mm}({\cal
P})=H^2_{mm}({\cal P})=0$, the map $\zf_{\sh}$ is an isomorphism
between $H^1_{mm}({\cal S})$ and $H^2_{mm}({\cal R})$. So, if
$a\neq 0,$ both spaces vanish. In the case $a=0,$ we obtain, see
Formulary 2, $H^2_{mm}({\cal R})=x^{m-2}\lp\R
x\p_{23}+\R(x\p_{31}+(m-1)y\p_{23})\rp$.  \rule{1.5mm}{2.5mm}\\

\noindent {\bf Remark.}  So, when passing from the ${\cal P}$- to
the ${\cal R}$-cohomology, the ${\cal P}$-class of degree
$(kr)=(11)$ is lost for $a\neq 0$, whereas for $a=0$ new $\cal
R$-classes of degree $(kr)=(mm)$, $m\ge 2$ appear. This is
coherent with the picture described in section \ref{seshds}, if
$\zf(\zs)$, $\zs=\frac{x^{m-1}}{D}(gx+hy)Y_3\in {\cal
S}^1_{mm}\cap ker\p_{\cal S},$ $m\ge 2,$ $g,h\in\R$, is a
coboundary in ${\cal P}$ but not in ${\cal R}$. Since $\pi_{\cal
S}\p_{\cal P}\zs=\p_{\cal S}\zs=0,$ we definitely have
$\zf(\zs)=\p_{\cal P}\zs$. If $\zf(\zs)$ were an ${\cal
R}$-coboundary, $\zf(\zs)=\p_{\cal R}\zr,$ $\zr\in {\cal
R}^1_{mm},$ the difference $\zr-\zs$ would be a ${\cal
P}$-coboundary in view of Theorem \ref{H1PR}: $\zr-\zs=\p_{\cal
P}p$, $p\in {\cal P}^0_{mm}=0$.


\section{3 - cohomology spaces}

\subsection{$\mathbf{{\cal S}}$ - cohomology}

\begin{prop}
The third cohomology group of ${\cal S}$ vanishes: $H^3({\cal
S})=0.$
\end{prop}
{\it Proof.} Obvious.  \rule{1.5mm}{2.5mm}

\subsection{$\mathbf{{\cal P}}$ - cohomology}

\begin{prop}(i) For $a\neq 0,$ \[\begin{array}{lll}H^3({\cal P})&=&H^3_{23}({\cal P})\oplus H^3_{00}({\cal
P})\\&=&\R Y_{123}\oplus \R\p_{123},\end{array}\](ii) for $a=0$,
\[\begin{array}{lll}H^3({\cal P})&=&\bigoplus_{k=1}^{\infty}H^3_{2k,3k}({\cal P})\oplus
\bigoplus_{m=0}^{\infty}H^3_{0m}({\cal
P})\\&=&\bigoplus_{m=0}^{\infty}D^{m}\R
Y_{123}\oplus\bigoplus_{m=0}^{\infty}\R
x_3^m\p_{123}.\end{array}\]
\end{prop}
{\it Proof.} Of course any cochain $\pi=qY_{123}\in {\cal
P}^3_{kr},$ $q\in Q_{kr},$ $k\ge 0,r\ge 0,k\le r$ is a cocycle.
For $k\ge 0,r\ge 1,k\le r$, this cocycle is a coboundary, if there
are $\zvr_1,\zvr_2,\zvr_3\in Q_{kr}$ such that
\[X_1(\zD_{k0}\zvr_1)+X_2(\zD_{k0}\zvr_2)+X_3(\zD_{kr}\zvr_3)=q.\]
\noindent {\bf 1. $\mathbf{r=0}$}\\

\noindent The cocycle reads
\[\pi^3_{00}=\frac{\za_0}{D}Y_{123}=\za_0\p_{123}\;\;\;(\za_0\in\R)\]
and no coboundaries do exist.\\

\noindent {\bf 2. $\mathbf{k=0,r\ge 1}$}\\

\noindent The coboundary condition reads $X_3(\zvr_3)=q$. If
$a\neq 0,$ we set $\zvr_3=X_3^{-1}(q).$ Otherwise we decompose $q$
in the form $q=\frac{\za_1z^r}{D}+q',$ $\za_1\in\R,q'\in im\,X_3.$
Then the cocycle
\[\pi^3_{0r}=\frac{\za_1z^r}{D}Y_{123}=\za_1z^r\p_{123}\;\;\;(a=0,r\ge 1,\za_1\in\R)\]
is not a coboundary (except of course if $\za_1=0$).\\

\noindent {\bf 3. $\mathbf{k\ge 1,r\ge 1,k=r}$}\\

\noindent Coboundary condition: $X_1(\zvr_1)+X_2(\zvr_2)=q$. As
$2r-3k=-r\neq 0,$ it suffices to take $\zvr_1=0$ and
$\zvr_2=X_2^{-1}(q)$.\\

\noindent {\bf 3. $\mathbf{k\ge 1,r\ge 1,k\le r-1}$}  (then $r\ge 2$)\\

\noindent Condition: $X_1(\zvr_1)+X_2(\zvr_2)+X_3(\zvr_3)=q$. \\

\noindent {\bf 3.a. $\mathbf{2r-3k\neq 0}$}\\

\noindent We only need choose $\zvr_1=\zvr_3=0$ and
$\zvr_2=X_2^{-1}(q).$\\

\noindent {\bf 3.b. $\mathbf{2r-3k=0}$}  (then $(kr)\in\{(23),(46),(69),\ldots\}$)\\

\noindent Here $X_2=0$ and $X_3=-2X_1$ and the condition reads
$X_1(\zvr_1-2\zvr_3)=q$. If $a\neq 0$ and $(kr)\neq (23)$, take
$\zvr_3=0$ and $\zvr_1=X_1^{-1}(q)$. If $a=0$ or $(kr)=(23)$, set
again $q=\za_2D^{\frac{k}{2}-1}+q'$, $\za_2\in\R,q'\in im\,X_1.$
The cocycle
\[\pi^3_{2m,3m}=D^{m-1}\za_2Y_{123}\;\;\;(a=0,m\in\N^*\mbox{ or }a\neq 0,m=1;\za_2\in\R)\]
can not be a coboundary (if $\za_2\neq 0$).

\subsection{$\mathbf{{\cal R}}$ - cohomology}

\begin{theo} (i) If $a\neq 0$, \[\begin{array}{lll}H^3({\cal R})&=&H^3_{23}({\cal R})\oplus H^3_{00}({\cal R})\\
&=&\R Y_{123}\oplus\R\p_{123},\end{array}\] (ii) if $a=0,$
\[\begin{array}{lll}H^3({\cal R})&=&\bigoplus_{k=1}^{\infty}H^3_{2k,3k}({\cal R})\oplus
\bigoplus_{m=0}^{\infty}H^3_{0m}({\cal
R})\oplus\bigoplus_{m=1}^{\infty}H^3_{mm}({\cal
R})\\&=&\bigoplus_{m=0}^{\infty}D^{m}\R
Y_{123}\oplus\bigoplus_{m=0}^{\infty}\R
x_3^m\p_{123}\oplus\bigoplus_{m=0}^{\infty}x_1^m(\R x_1+\R x_2
)\p_{123}.\end{array}\]
\end{theo}


\noindent {\it Proof.}  If $a\neq 0,(kr)\neq (11)$ and if
$a=0,(kr)\neq (mm),m\ge 1$, then \[i_{\sh}\in Isom(H^3_{kr}({\cal
R}),H^3_{kr}({\cal P}))\;\;\;(a\neq 0,(kr)\neq (11)\mbox{ or
}a=0,(kr)\neq (mm),m\ge 1).\]

As $H^3_{mm}({\cal P})=0,$ $m\ge 1,$ it follows from Corollary
\ref{HS0} that $\zf_{\sh}\in Isom(H^2_{mm}({\cal
S})/ker^2_{mm}\zf_{\sh},H^3_{mm}({\cal R}))$. If $a=0,m\ge 2,$ the
group $H^2_{mm}({\cal P})$ also vanishes and \[\zf_{\sh}\in
Isom(H^2_{mm}({\cal S}),H^3_{mm}({\cal R}))\;\;\;(a=0,m\ge 2).\]
For $m=1$, we have to compute the kernel $ker^2_{mm}\zf_{\sh}$. If
$\zf_{\sh}[\zs]_{\cal S}=0,$ $\zs\in {\cal S}^2_{11}\cap
ker\,\p_{\cal S}={\cal S}^2_{11}$, the image $\zf(\zs)=\p_{\cal
P}\zs$ is an ${\cal R}$-coboundary $\p_{\cal R}\zr$, $\zr\in {\cal
R}^2_{11}.$ So, $\zr-\zs\in {\cal P}^2_{11}\cap ker\,\p_{\cal P}.$
Equation \ref{pi211} shows that
\[\zr-\zs=\zp^2_{11}=\frac{1}{D}\left[\lp(\za-\frac{a}{b}\zb) x+\zb y\rp Y_{23}+
\lp\zb x+(\frac{a}{b}\zb-\za)y\rp
Y_{31}\right]+\frac{x}{D}\left[\frac{a}{b}\zb Y_{23}+
\frac{a}{b}\za Y_{31}\right],\] where the first (resp. second)
term of the r.h.s. is a member of ${\cal R}^2_{11}$ (resp. ${\cal
S}^2_{11}$). Hence, if $a=0$, the cochain $\zs$ vanishes and
$ker^2_{11}\,\zf_{\sh}=0$: \[\zf_{\sh}\in Isom(H^2_{11}({\cal
S}),H^3_{11}({\cal R}))\;\;\;(a=0).\] If $a\neq 0$ and
$\zs=\frac{x}{D}(iY_{23}+jY_{31})$, see Formulary 2, we get
$\za=-\frac{b}{a}j,\zb=-\frac{b}{a}i$, which provides $\zr$ in
terms of $\zs$. Equations \ref{X_1}, \ref{X_2}, and \ref{f} then
immediately give $ker^2_{11}\,\zf_{\sh}=H^2_{11}({\cal S})$, so
that \[H^3_{11}({\cal R})=0\;\;\;(a\neq 0).\] The announced
theorem is then a direct consequence of Equation \ref{f}.
\rule{1.5mm}{2.5mm}

\section{Further results\label{coho7}}

In this section, we provide complete results regarding the formal
cohomology of structure 7 of the DHC,
\[\zL_7=b(x_1^2+x_2^2)\p_1\w\p_2+((2b+c)x_1-ax_2)x_3\p_2\w\p_3+(ax_1+(2b+c)x_2)x_3\p_3\w\p_1.\]
We obtained these upshots, using the same method as above for
structure 2. Computations are quite long and will not be published
here. We assume that $c\neq 0$, otherwise we recover structure 2.

In the following theorems, the $Y_i$ ($i\in\{1,2,3\}$) denote the
same vector fields as above, namely,
\[Y_1=x_1\p_1+x_2\p_2,Y_2=x_1\p_2-x_2\p_1,Y_3=x_3\p_3.\]
Moreover, we set \[D'=x_1^2+x_2^2, D=(x_1^2+x_2^2)x_3.\] If
$\frac{b}{c}\in\Q, b(2b+c)<0,$ we denote by $(\zb,\zg)\simeq
(b,c)$ the irreducible representative of the rational number
$\frac{b}{c}$, with positive denominator, $\zb\in\Z,\zg\in\N^*.$
If $\frac{b}{c}\in\Q, b(2b+c)>0,$ $(\zb,\zg)\simeq (b,c)$ denotes
the irreducible representative with positive numerator,
$\zb\in\N^*,\zg\in\Z^*.$ Furthermore, we write $\zL$ instead of
$\zL_7$, $\bigoplus_{ij}\op{Cas}(\zL) Y_{ij}$ instead of
$\op{Cas}(\zL) Y_{23}+\op{Cas}(\zL) Y_{31}+\op{Cas}(\zL) Y_{12}$,
$\op{Sing}(\zL)=\bigoplus_{r\ge 0}\op{Sing}^r(\zL)$ instead of
$\R[[x_3]]=\bigoplus_{r\ge 0}\R x_3^r$, and $C_{\zg}Y_3$
($\zg\in\{2,4,6,\ldots\}$) instead of $\R
D'^{\frac{\zg}{2}-1}x_3^{-1}Y_3=\R D'^{\frac{\zg}{2}-1}\p_3.$

Some comments on the results given in the theorems hereafter can
be found below.

\begin{theo} If $a\neq 0$, the cohomology spaces are $$H^0(\zL)=\op{Cas}(\zL)=\R,\quad H^1(\zL)=\bigoplus_i\op{Cas}(\zL) Y_i,$$ $$H^2(\zL)=
\bigoplus_{ij}\op{Cas}(\zL)Y_{ij},\quad H^3(\zL)=\op{Cas}(\zL)
Y_{123}\oplus \op{Sing}^0(\zL)\p_{123}$$ \end{theo}

\begin{theo} If $a=0$ and $b=0$, the cohomology is $$ H^0(\zL)=\op{Cas}(\zL)=\bigoplus_{r\ge 0}\R D'^r, \quad
H^1(\zL)=\bigoplus_i\op{Cas}(\zL)Y_i,$$
$$ H^2(\zL)= \bigoplus_{ij}\op{Cas}(\zL)Y_{ij}\oplus\op{Sing}(\zL)\p_{12}, \quad
H^3(\zL)=
\op{Cas}(\zL)Y_{123}\oplus\op{Sing}(\zL)\p_{123}$$\end{theo}

\begin{theo} If $a=0$ and $2b+c=0$, the cohomology groups are $$ H^0(\zL)=\op{Cas}(\zL)=\bigoplus_{r\ge 0}\R
x_3^r,\quad H^1(\zL)=\bigoplus_i\op{Cas}(\zL) Y_i\oplus C_2Y_3,$$
$$ H^2(\zL)=\bigoplus_{ij}\op{Cas}(\zL)Y_{ij}\oplus\op{Sing}(\zL)\p_{12}\oplus C_2Y_3\w (\R Y_1+\R Y_2),$$
$$ H^3(\zL)=\op{Cas}(\zL)Y_{123}\oplus\op{Sing}(\zL)\p_{123}\oplus C_2Y_3\w \R Y_{12}$$\end{theo}

\begin{theo} If $a=0$ and $\frac{b}{c}\notin\Q$ or $\frac{b}{c}\in\Q,b(2b+c)<0$,
$$ H^0(\zL)=\op{Cas}(\zL)=\R,\quad H^1(\zL)= \bigoplus_i\op{Cas}(\zL)Y_i\;\oplus\begin{cases}(b,c)\simeq
(-1,\zg),\zg\in\{4,6,8,\ldots\}:C_{\zg}Y_{3}\\otherwise:0\end{cases},$$
$$H^2(\zL)=\bigoplus_{ij}\op{Cas}(\zL)Y_{ij}\oplus \op{Sing}(\zL)\p_{12}\oplus \begin{cases}(b,c)\simeq
(-1,\zg),\zg\in\{4,6,8,\ldots\}:\\C_{\zg}Y_3\w(\R Y_1+\R
Y_2)\\otherwise:\\0\end{cases},$$
$$ H^3(\zL)=\op{Cas}(\zL)Y_{123}\oplus\op{Sing}(\zL)\p_{123}\oplus\begin{cases}(b,c)\simeq(-1,\zg),\zg\in\{4,6,8,\ldots\}:C_{\zg}Y_3\w\R Y_{12}\\\mbox{otherwise}:0\end{cases}$$
\end{theo}

\begin{theo} \label{a0p>0 } If $a=0$ and $\frac{b}{c}\in\Q,b(2b+c)>0$,
\[H^0(\zL)=\op{Cas}(\zL)=\bigoplus_{n\in\N,n\zg\in
2\Z}\R D'^{n\zb+\frac{n\zg}{2}}x_3^{n\zb},\quad
H^1(\zL)=\bigoplus_i\op{Cas}(\zL)Y_i,\]
$$H^2(\zL)=\bigoplus_{ij}\op{Cas}(\zL)Y_{ij}\oplus\op{Sing}(\zL)\p_{12},\quad H^3(\zL)=\op{Cas}(\zL)Y_{123}\oplus\op{Sing}(\zL)\p_{123}$$
\end{theo}

{\bf Remark}. The preceding results allow to ascertain that
Casimir functions are closely related with Koszul-exactness or
``quasi-exactness'' of the considered structure. Observe that
$C_{\zg}=\R D'^{-1+\frac{\zg}{2}} x_3^{-1}$
($\zg\in\{2,4,6,\ldots\}$) has the same form as the basic Casimir
in Theorem \ref{a0p>0 } and that the negative superscript in
$x_3^{-1}$ can only be compensated via multiplication by
$Y_3=x_3\p_3$. Clearly, such a compensation is not possible for
$\R D'^{-n+\frac{n\zg}{2}} x_3^{-n}$, $n>1$. Hence cocycle
$C_{\zg}Y_3$ is in some sense ``Casimir-like'' and ``accidental''.
Note eventually that the ``weight'' of the singularities in
cohomology increases with closeness of the considered Poisson
structure to Koszul-exactness.

\section{Suitable family of quadratic structures}

The objective of this final section is to explain that our
technique applies to all the quadratic Poisson classes induced by
a special type of $r$-matrices.\\

It is well-known that the action tangent to the canonical action
of the Lie group $GL(n,\R)$ on $\R^n$, is the Lie algebra
homomorphism \[J:g:=gl(n,\R)\ni a=(a_{ij})\raa a_{ij}x_i\p_j\in
Sec(T\R^n)=:{\cal X}^1(\R^n)\] ($Sec(T\R^n)$: Lie algebra of
smooth sections of the tangent bundle $T\R^n$, i.e. Lie algebra of
vector fields of $\R^n$; $x=(x_1,\ldots,x_n)$: canonical
coordinates of $\R^n$; $\p_1,\ldots,\p_n$: partial derivatives
with respect to these coordinates), which is a Lie isomorphism if
valued in the Lie algebra ${\cal X}^1_0(\R^n)$ of linear vector
fields.

Let us recall that a standard construction allows to associate to
any Lie algebroid $E$ a Gerstenhaber algebra, made up by a graded
Poisson-Lie algebra structure on the shifted Grassmann algebra
$Sec(\bigwedge E)[1]$ of multi-sections of $E$. We denote this
Schouten-Nijenhuis superbracket, which extends the algebroid
bracket on $Sec(E)$, by $[.,.]_E$ or simply by $[.,.]$, if no
confusion is possible.

The above Lie homomorphism $J$ extends to a Gerstenhaber
homomorphism
\[J:\bigwedge g\raa {\cal X}(\R^n):=Sec(\bigwedge T\R^n),\]
where the Gerstenhaber structures have been obtained as just
mentioned. Let \[\tilde{\bigwedge}\R^n=\bigoplus_k\lp {\cal
S}^k(\R^n)^*\otimes\bigwedge\,^k\R^n\rp\] be the Gerstenhaber
subalgebra---of the algebra ${\cal X}(\R^n)$ of poly-vector
fields---made up by $k$-vectors with coefficients in the
corresponding space of homogeneous polynomials in $x\in\R^n$. It
is obvious that $J$ viewed as Gerstenhaber homomorphism with
target algebra $\tilde\bigwedge\R^n$, \[J:\bigwedge
g\raa\tilde{\bigwedge}\R^n,\] is onto. It is also known that the
restriction $J^k$ to the space $\bigwedge^kg$ has a non-trivial
kernel (provided that $k,n\ge 2$).\\

Remember that a classical $r$-matrix is a bi-matrix $r\in g\wedge
g$ that verifies the classical Yang-Baxter equation $[r,r]=0$. The
space ${\cal S}^2(\R^n)^*\otimes\bigwedge^2\R^n$ of quadratic
bi-vectors coincides with the image $J^2(g\wedge g)$ and
\[J^3[r,r]=[J^2r,J^2r],\;r\in g\wedge g.\] So any bi-vector
$\zL=J^2r$ that is the image of an $r$-matrix is a quadratic
Poisson structure of $\R^n$. Conversely, any quadratic Poisson
tensor of $\R^n$ is induced by at least one bi-matrix. However,
the characterization of those quadratic Poisson structures that
are implemented by an $r$-matrix, is an open problem (see
\cite{MMR}). Quadratic Poisson structures generated by an
$r$-matrix are of importance e.g. in deformation
quantization, in particular in view of Drinfeld's method.\\

One easily understands that the study of $r$-matrix induced
structures involves the orbit $O_{\zL}$ and stabilizer $G_{\zL}$
of the considered structure $\zL$ for the canonical action on
tensors of the general linear group $GL(n,\R)$. Our paper is
confined to the three-dimensional Euclidean setting. Let us recall
that the isotropy group $G_{\zL}$ is a Lie subgroup of
$G:=GL(3,\R)$, the Lie algebra of which is the stabilizer
$g_{\zL}$ of $\zL$ for the corresponding infinitesimal action,
\[g_{\zL}=\{a\in gl(3,\R): [Ja,\zL]=0\}.\]

As already stated, the objective of this paper is to provide a
universal approach to the formal Poisson cohomology for a broad
family of quadratic structures in the classical physical space.
Let us just indicate that the classes of the DHC that are
accessible to our modus operandi are exactly those classes that
are implemented by $r$-matrices in $g_{\zL}\wedge g_{\zL}$. The
quadratic Poisson tensors can defacto be classified according to
their membership in the family of structures induced by an
$r$-matrix in the ``stabilizer''. It then turns out that the
members of this family are those tensors that read as a linear
combination of wedge products of mutually commuting linear vector
fields, hence those that are accessible to the above detailed and
applied technique.

We refer to the classes implemented by an $r$-matrix in the
stabilizer as admissible classes. The DHC gives the quadratic
Poisson structures up to linear transformations. Let us finally
mention that admissibility is of course effectively independent of
the chosen representative. This means that if $\zL=J^2r,r\in
g_{\zL}\w g_{\zL},[r,r]=0,$ then any equivalent structure
$A_*\zL,A\in G$---$A_*$ denotes the natural action of $A$---has
the same property, i.e. $A_*\zL$ is also induced by an $r$-matrix
and this matrix can be chosen in the stabilizer of $A_*\zL$.
Indeed, it is easily seen that the orbit of $J^2r$ for the
canonical $G$-action, is nothing but the pointwise image of the
orbit of $r$ for the adjoint action $\op{Ad}$ of $G$,
\[A_*(J^2r)=J^2(\op{Ad}(A)r).\] Since the
adjoint action respects the Schouten-Nijenhuis bracket,
$\op{Ad}(A)r$ is an $r$-matrix in the stabilizer of $A_*\zL.$

\section{Acknowledgements}

M. Masmoudi is grateful to A. Roux for fruitful discussions at the
Paul Verlaine University of Metz.

\end{document}